\documentclass[a4paper,12pt]{article}
\usepackage[warn]{mathtext}
\usepackage{cmap}
\usepackage[T2A]{fontenc}
\usepackage[cp1251]{inputenc}
\usepackage[russian]{babel}
\usepackage{amsmath}
\usepackage{wasysym}
\usepackage{amssymb,amsthm}
\usepackage[dvips]{graphicx}

\theoremstyle{plain}
\newtheorem{theorem}{Теорема}[section]
\newtheorem{corollary}[theorem]{Следствие}
\newtheorem{lemma}[theorem]{Лемма}
\newtheorem{proposition}[theorem]{Предложение}

\theoremstyle{definition}
\newtheorem{definition}[theorem]{Определение}

\def\Evol{\operatorname{Evol}}

\begin{document}

\title{
Периодичность схем Рози для подстановочных слов}

\author{А.~Я.~Белов\footnote{ работа выполнена при частичной финансовой поддержке гранта РФФИ № 14-01-00548}, И.~Митрофанов\footnote{ работа выполнена при частичной финансовой поддержке фонда Дмитрия Зимина ``Династия" и гранта РФФИ № 14-01-00548}}
\date{}
\maketitle

\section{Введение.}

Пусть $M$ --- компакт, $U\subset M$
--- его открытое подмножество,
 $f\colon M \to M$ --- гомеоморфизм компакта в себя и $x_0\in M$ --- начальная точка.
 Положим $ w_n= a,\text{если $f^n(x_0)\in U$}$, $w_n=
b,\text{если $f^n(x_0)\notin U$}$;  Слово $w$ называется {\em
эволюцией} точки $x$. Символическая динамика исследует взаимосвязь
свойств динамической системы $(M,f)$ и комбинаторных свойств слова
$w$. Для слов над алфавитом, состоящим из большего числа символов,
нужно рассмотреть несколько характеристических множеств. Для
уточнения абстрактной постановки задачи исследования соответствия
комбинаторных  и топологических свойств слов очень важна модельная
ситуация, от которой можно отталкиваться при дальнейших
исследованиях. Проблематика, связанная с построением слов Штурма,
очень важна в комбинаторике слов и в алгебре. Известна
классическая

\begin{theorem}[Теорема эквивалентности]
Пусть $W$ -- бесконечное рекуррентное слово над бинарным алфавитом
$A=\{a,b\}$. Следующие условия  эквивалентны с точностью до
счетного множества последоательностей: (1) Слово $W$ является
словом Штурма, то есть количество различных подслов длины $n$
слова $W$ равно $T_n(W)=n+1$ для любого $n\geq 1$. (2) Слово не
периодично и является {\em сбалансированным}, то есть для любых
двух подслов $u,v\subset W$ одинаковой длины выполняется
неравенство $||v|_a-|u|_a|\leq 1$, где $|w|_a$ обозначает
количество вхождений символа $a$ в слово $w$. (3) Слово $W=(w_n)$
является {\em механическим словом} с иррациональным $\alpha$, то
есть существуют такое иррациональное $\alpha$, $x_0 \in [0,1]$ и
интервал $U\subset \mathbb{S}^1$, $|U|=\alpha$ такие, что $ w_n=a$
если ${T_{\alpha}}^n(x_0)\in U$, иначе $w_n=b$. (4) Слово $W$
получается путем предельного перехода последовательности слов,
каждое из которых получается из предыдущего путем подстановки вида
$a^kb\to b, a^{k+1}b\to a$ либо подстановки вида $b^ka\to a,
b^{k+1}a\to b$. Показатель $k$ зависит от шага и отвечает остатку
разложения $\alpha$ в цепную дробь. Если эти показатели $k_i$
периодически повторяются, то $\alpha$ есть квадратичная
иррациональность. (5) Сверхслово $W$ р.р. и имеет
последовательность графов Рози с одной входящей и одной исходящей
развилкой.
\end{theorem}

На основе каждого из свойств (1)--(5) разными авторами строятся
обобщения последовательностей Штурма. Подробнее -- см.
\cite{BelKondMitr}. Там же обсуждаются персоналии. Рассмотрение
общего случая слов, у которых графы Рози имеют большее число
развилок (но все же их число конечно), т.е. слов с линейной
функцией сложности приводит к изучению слов, порождаемых {\bf
перекладыванием отрезков.} Перекладывания отрезков, т.е. кусочно
непрерывные преобразования одномерного комплекса, естественным
образом служат обобщением вращения круга. (Сдвиг окружности
эквивалентен  перекладыванию двух отрезков с сохранением
ориентации.)

Имеется комбинаторный критерий (для общего случая, не обязательно
являющегося регулярным) того, что данное сверхслово является
перекладыванием отрезков \cite{BelovChernInt}. Критерий таков:
ребра развилок в графах Рози помечены буквами ``п'' и ``л'',
отвечающими подразбиению соответствующих отрезков на правые и левые
части. При возникновении биспециального слова (двусторонней
развилки) одна из возможностей ``пл'' либо ``лп'' должна быть
уничтожена. В процессе эволюции стрелка наследует ту же букву.

\begin{theorem}[\cite{BelovChernInt}]    \label{ThBelCher}
Равномерно-рекуррентное слово $W$
 порождается перекладыванием отрезков тогда и только тогда, когда
слово обеспечивается асимптотически правильной эволюцией
размеченных графов Рози.
\end{theorem}

Аналогичное утверждение можно сформулировать для перекладывания
отрезков с изменением ориентации. В работе \cite{FZ3} был
независимо  получен другой критерий порождаемости слов
преобразованием перекладывания отрезков с сохранением ориентации,
удовлетворяющих следующему условию: траектория каждой концевой
точки отрезка перекладывания не попадает на концевую отрезка
перекладывания, в том числе сама на себя.

Графы Рози представляются комбинаторным языком, наиболее
приспособленным для задания одномерных динамических систем. Ещё
одним комбинаторным языком является язык подстановок.
Подстановочные динамические системы и тесно связанные с ними
$D0L$-системы интенсивно исследовались рядом авторов. В этой связи
следует обратиться к характеризации слов Штурма через свойство
(4). Что можно сказать про динамические системы для произвольной
системы подстановок? Для поворота окружности периодичность может
наблюдаться только для случая ранга $2$, т.е. для подстановок,
собственные значения соответствующей матрицы которых есть
квадратичная иррациональность. Для перекладывания отрезков в общем
случае возможно появление высших иррациональностей. В этой связи
важно уметь переводить свойство периодичности на язык графов Рози.
Эту возможность предоставляет теорема, доказываемая в настоящей
статье:

\begin{theorem} \label{ThMain}
Непериодичное равномерно-рекуррентное слово является
подстановочным тогда и только тогда, когда протокол
детерминированной эволюции его схем Рози периодичен, возможно, с
предпериодом.
\end{theorem}

Обратимся к условиям (4) и (5) для графов Рози слов Штурма. Для
случая одной входящей (соответственно, одной исходящей) развилки
периодичность событий в схеме Рози означает периодичность
разложения в цепную дробь числа $\alpha$.
 Поэтому теорему  \ref{ThMain} можно
рассматривать как обобщение теоремы Лагранжа о периодичности
разложения квадратичной иррациональности в цепную дробь на высшие
иррациональности. В случае размеченных схем Рози для
перекладывания отрезков периодичность также равносильна
подстановочности.

Из теоремы \ref{ThMain} также вытекает обобщение классической
теоремы А.А.Маркова. Имеется счетная последовательность
алгебраических констант $\{\alpha_i\}\to 1/3$ и чисел
$\{\beta_i\}$ такие, что для любого алгебраического числа
$\gamma$, не эквивалентного $\beta_j$ для всех $j<i$ неравенство
$|\gamma-p/q|<1/cq^2$ имеет конечное число решений для всех
$c<\alpha_i$. Величина наилучшего приближения определяется
отношением минимального интервала к максимальному в процессе
алгоритма Евклида на паре отрезков. Индукция Рози обобщает
Алгоритм Евклида. Для перекладывания $\le 4$ отрезков
соответствующие результаты были получены в работе \cite{Ferenci2}.
Отношение длин отрезков равно отношению частот соответствующих
подслов. Назовем такое отношение {\em дискретной частью спектра},
если такое же или меньшие отношения достигаются на счетном
множестве эволюций схем Рози. Аналогично определяется {\em
дискретная часть спектра} для перекладывания отрезков.

\begin{corollary}[Теорема Маркова для произвольного числа развилок]
Дискретная часть спектра достигается на периодических схемах Рози.
Соответствующие отношения суть алгебраические числа. Аналогичный
факт верен для размеченных схем Рози реализующий перекладвание
отрезков.
\end{corollary}

Перекладывания отрезков возникают при изучении потоков на
поверхностях отрицательной кривизны, в частности, из изучения
бильярда с рациональными углами. В этом случае также легко
сооружается подстановочная система. Если осуществлять эволюцию
случайным образом, всякий раз выбирая возможность с вероятностью
$1/2$,  то ситуации для слов Штурма отвечает инвариантная
относительно перехода к следующему остатку Гауссова мера, а для
перекладывания отрезков -- инвариантная относительно индукции Рози
мера на пространствах Тейхмюллера.

Из доказательства теоремы \ref{ThMain} вытекает также решение
известных открытых вопросов \cite{SemenovReview}.

\begin{theorem}
Существует алгоритм проверки  равномерно-рекуррентности
морфического слова $h(\varphi^\infty(s))$.
\end{theorem}

\begin{theorem}
Существует алгоритм проверки  периодичности
морфического слова $h(\varphi^\infty(s))$.
\end{theorem}

Из теоремы \ref{ThMain} вытекает теорема Вершика-Лившица о
периодичности диаграмм Брателли для марковских компактов,
порожденных подстановочными системами \cite{VL}. Значительно более
сложной, чем в теореме Вершика-Лившица, частью рассуждений
является доказательство периодичности схем Рози для морфических
слов. Подстановка может быть плохо устроена: образы различных букв
могут содержать друг друга. Для некоторых подстановочных систем
схемы Рози отвечают графам Рози, в которых простые пути между
развилками заменяются на ориентированные рёбра длины $1$. В
частности, такова ситуация в случае {\em равноблочных
маркированных циркулярных} $DOL$-последовательностей, изученном в
\cite{Frid}. В общем случае {\em схема Рози} являет собой граф,
различные части которого взяты из графов Рози разных порядков, на
рёбрах графа по некоторому правилу написаны слова. Назовем  {\em
весом ребра} в схеме Рози длину соответствующего слова. Ключевым
местом является доказательство того, что {\it если слово $W$
порождено примитивным морфизмом, то отношения весов во всех схемах
Рози ограничены}. (То же верно для р.р. морфических слов.) Далее
можно воспользоваться результатом J.Cassaigne \cite{Cassaigne} о
том что если $W$ равномерно рекуррентно и $\liminf
T_W(n)/n<\infty$, то $\limsup T_W(n+1)-T_W(n)<\infty$ и установить
ограниченность числа вершин в схемах Рози, получающихся путем
элементарной последовательности эволюций из заданной.

\paragraph{\bf Основные определения.}
  Если на $i$-том месте в слове $u$ стоит буква $a$, то
положим $u[i]=a$, $|u|$ есть  {\em длина} слова $u$. На словах
существует естественная структура частично упорядоченного
множества: $u_1\sqsubseteq u_2$, если $u_1$ является подсловом
$u_2$. Будем обозначать $u_1\sqsubseteq_k u_2$, если слово $u_1$
входит в $u_2$ хотя бы $k$ раз. Сверхслово $W$ называется {\em
рекуррентным}, если любое его подслово встречается в $W$
бесконечно много раз и {\em равномерно рекуррентным}, если для
любого его подслова $v$ существует такое число $k(v,W)$, что если
$u\sqsubseteq W$ и $|u|\geq k(v,W)$, то $v\sqsubseteq u$.
Множество слов $A^{*}$ над алфавитом $A$ является свободным
моноидом с операцией конкатенации и единицей -- пустым словом.
Отображение $\varphi:A^*\to B^*$ называют {\em морфизмом}, если
оно сохраняет структуру моноида. Если алфавиты $A$ и $B$ совпадают
и существует такая буква $a_1$, что $\varphi(a_1)=a_1u$ для
некоторого слова $u$ и $\varphi^k(u)$ не является пустым словом ни
для какого $k$, то бесконечное слово
$a_1u\varphi(u)\varphi^2(u)\varphi^3(u)\varphi^4(u)\dots $
называется {\em чисто морфическим}, пишут
$W=\varphi^{\infty}(a_1)$. Если существует такая степень морфизма
$\varphi^k$, что для любых двух букв $a_i$ содержится в
$\varphi^k(a_j)$, морфизм называют {\em примитивным}.

\section{Графы Рози, графы со словами и схемы Рози}

 Для бесконечного вправо слова $W$ определёны {\em графы Рози}.
Граф Рози {\em порядка $k$} обозначается $G_k(W)$, (или просто
$G_k$). Вершины его соответствуют всевозможным различным подсловам
длины $k$ сверхслова $W$. Две вершины графа $u_1$ и $u_2$
соединяются направленным ребром, если в $W$ есть такое подслово
$v$, что $|v|=k+1$, $v[1]v[2]\dots v[k]=u_1$ и $v[2]v[3]\dots
v[k+1]=u_2$. Если $w$ -- подслово $W$ длины $k+l$, ему
соответствует в $G_k$ путь длины $l$, проходящий по рёбрам,
соответствующим подсловам слова $w$ длины $k+1$.

{\em Графом со словами} будем называть связный ориентированный
граф, у которого на каждом ребре которого написано по два слова --
{\em переднее} и {\em заднее}, а кроме того, каждая вершина либо
имеет входящую степень $1$, а исходящую больше $1$, либо входящую
степень больше $1$ и исходящую степень $1$. Вершины первого типа
назовём {\em раздающими}, а второго -- {\em собирающими}.  {\em
Симметричный путь в графе со словами} -- это путь, первое ребро
которого начинается в собирающей вершине, а последнее ребро
кончается в раздающей. Каждый путь можно записать словом над
алфавитом, являющимся множеством рёбер графа. Оно называется {\em
рёберной записью пути}. Ребро пути $s$, идущее $i$-тым по счёту,
будем обозначать $s[i]$. Естественно определены отношения {\em
подпути} (пишем $s_1\sqsubseteq s_2$), {\em начала} и {\em конца}.
$s_1\sqsubseteq_k s_2$, если для соответствующих слов $u_1$ и
$u_2$ -- рёберных записей путей $s_1$ и $s_2$ -- выполнено
$u_1\sqsubseteq_k u_2$. Если последнее ребро пути $s_1$ идёт в ту
же вершину, из которой выходит первое ребро пути $s_2$, путь,
рёберная запись которого является конкатенацией рёберных записей
путей $s_1$ и $s_2$, будем обозначать $s_1s_2$. Определения пути,
подпути, рёберной записи, начала и конца имеют смысл для любых
ориентированных графов.

Введём понятие {\em переднего слова $F(s)$, соответствующего пути
$s$ в графе со словами.} Пусть $v_1v_2\dots v_n$ -- рёберная
запись пути $s$. В $v_1v_2\dots v_n$ возьмём
подпоследовательность: включим в неё $v_1$, а также те и только те
рёбра, которые выходят из раздающих вершин графа. Эти рёбра
назовём {\em передними образующими для пути $s$}. Возьмём передние
слова этих рёбер и запишем их последовательную конкатенацию, так
получаем $F(s)$. Аналогично определяется понятие {\em заднего
слова $B(s)$}. В пути $v_1v_2\dots v_n$ возьмём рёбра, входящие в
собирающие вершины и ребро $v_n$ в порядке следования -- это {\em
задние образующие для пути $s$}. Тогда последовательной
конкатенацией задних слов этих рёбер получается {\em заднее слово
$B(s)$ пути $s$}.

\begin{definition} \label{Def1}
Граф со словами будет являться {\em схемой Рози} для рекуррентного
непериодичного сверхслова $W$, если он удовлетворяет следующим
свойствам:

\begin{enumerate}
    \item Граф сильносвязен и состоит более чем из одного ребра.
    \item Все рёбра, исходящие из одной раздающей вершины графа, имеют
передние слова с попарно разными первыми буквами.  Все рёбра,
входящие в одну собирающую вершину графа, имеют задние слова с
попарно разными последними буквами.

    \item Для любого симметричного пути, его переднее и заднее слова
совпадают.
    \item Если есть два симметричных пути $s_1$ и $s_2$ и
выполнено $F(s_1)\sqsubseteq _k F(s_2)$, то $s_1\sqsubseteq _k
s_2$.
    \item Все слова, написанные на рёбрах графа, являются подсловами $W$.
     \item Для любого $u\sqsubset W$ существует симметричный путь,
слово которого содержит $u$.
     \item Для любого ребра $s$
существует такое слово $u_s$, принадлежащее $W$, что любой
симметричный путь, слово которого содержит $u_s$, проходит по
ребру $s$.
\end{enumerate}
Симмметричный путь $s$ называется {\em допустимым}, если
$F(s)\sqsubseteq W$.
\end{definition}

Пусть $W$ -- рекуррентное непериодичное сверхслово. Фиксируем
натуральное число $k$. Для упрощения дальнейшего считаем, что в
$W$ нет биспециальных слов длины ровно $k$. Рассмотрим $G_k$ --
граф Рози порядка $k$ для сверхслова $W$. Тогда  $G_k$ есть
сильносвязный орграф, не являющийся циклом. Построим граф со
словами $S$, вершинами которого будут специальные вершины графа
$G_k$, а рёбра -- простыми цепями, соединяющие специальные вершины
$G_k$. В этом случае {\it пути в графе $S$ соответствуют путям в
$G_k$, начинающимся и кончающимся в специальных вершинах. } Пусть
простой путь проходит в графе $G_k$ по вершинам
$$
a_{i_1}a_{i_2}\dots a_{k},\: a_{i_2}a_{i_3}\dots a_{k+1},\: \ldots
, a_{i_l}a_{i_{l+1}}\dots a_{i_{l+n-1}},\: \ldots,
\:a_{i_{n-k+1}}a_{i_{n-k}}\dots a_{i_n}
$$
и соединяет специальные вершины $a_{i_1}a_{i_2}\dots a_{k}$ и
$a_{i_{n-k+1}}a_{i_{n-k}}\dots a_{i_n}$. Сопоставим этому пути
{\em переднее} и {\em заднее} слова по следующему правилу:

\begin{enumerate}
\item Если вершина, соответствующая $a_{i_1}a_{i_2}\dots a_{i_k}$,
является раздающей, то {\em переднее} слово пути -- это
$a_{i_{k+1}}a_{i_{k+2}}\dots a_{i_n}$. \item Если вершина,
соответствующая $a_{i_1}a_{i_2}\dots a_{i_k}$, является в
собирающей, то {\em переднее} слово пути -- это
$a_{i_1}a_{i_2}\dots a_{i_n}$. \item Если вершина, соответствующая
$a_{i_{n-k+1}}a_{i_{n-k}}\dots a_{i_n}$, является собирающей, то
{\em заднее} слово пути -- это $a_{i_1}a_{i_2}\dots a_{i_n-k}$.
\item Если вершина, соответствующая $a_{i_{n-k+1}}a_{i_{n-k}}\dots
a_{i_n}$, является раздающей, то {\em заднее} слово пути -- это
$a_{i_1}a_{i_2}\dots a_{i_n}$.
\end{enumerate}

\begin{definition}
Если $S$ -- сильносвязный граф, не являющийся циклом, и $s$ --
путь в графе, то {\em естественное продолжение} пути $s$ {\em
вправо} -- это минимальный путь, началом которого является $s$ и
который оканчивается в раздающей вершине. {\em Естественное
продолжение} пути $s$ {\em влево} -- это минимальный путь, концом
которого является $s$ и который начинается в собирающей вершине.
\end{definition}

Для сильносвязных нецикличных графов естественное продолжение
существует всегда и единственно. Напишем слова на рёбрах $S$: для
каждого ребра в качестве переднего слова берётся переднее слово
того пути в $G_k$, который соответствует естественному расширению
вправо этого ребра. Заднее же слово -- это заднее слово пути в
$G_k$, соответствующего естественному продолжению влево
рассматриваемого ребра.

\begin{proposition} \label{Prop_1}
а) В полученном графе со словами $S$ для любого пути $s$ слово
$F(s)$ -- это переднее слово того пути графа $G_k$, который
соответствует естественному продолжению вправо пути $s$.
Аналогично, $B(s)$ -- это заднее слово того пути в $G_k$, который
соответствует естественному продолжннию пути $s$ влево (в графе
$S$). б) Пусть $S$ -- определённый выше граф со словами. Тогда он
является схемой Рози для сверхслова $W$.
\end{proposition}

Пусть $W$ -- рекуррентное непериодичное сверхслово, $S$ -- схема
Рози для этого сверхслова. Из сильносвязности $S$ следует, что
существует ребро $v$, идущее из собирающей вершины в раздающую.
Такие рёбра будем называть {\em опорными}. Естественным образом
определяется $S''$ -- {\em эволюция $(W,S,v)$ схемы Рози $S$ по
опорному ребру $v$}: само ребро $v$ уничтожается и заменяется на
минимальные пути через $v$ строго содержащие $v$, отвечающие
подсловам $W$. При этом естественным образом определяются передние
и задние слова для нового графа со словами. Построенный таким
образом граф со словами $S''$ назовём {\em элементарной эволюцией}
$(W,S,v)$.

Имеет место следующая

\begin{theorem} \label{T1}
Пусть $W$ -- схема Рози. Тогда элементарная эволюция $(W,S,v)$
является схемой Рози для сверхслова $W$ (то есть удовлетворяет
свойствам (1)--(7) определения \ref{Def1}).
\end{theorem}

К полученному графу со словами можно, выбрав некоторое ребро,
снова применить элементарную эволюцию, получая таким образом
последовательность схем Рози. На каждом шаге опорное ребро может
быть выбрано, вообще говоря, несколькими способами. Пронумеруем
рёбра исходной схемы и зафиксируем {\it метод эволюции} --
алгоритм, который по пронумерованной схеме определяет, какое
опорное ребро использовать, а после перенумеровывает рёбра в
проэволюционировавшей схеме. Тогда последовательность схем Рози
определяется однозначно. Потребуем также, чтобы метод эволюции не
зависел от слов, написанных на рёбрах, то есть работал с {\it
облегчённой нумерованной схемой}.

\begin{definition}
{\it Облегчённая нумерованная схема для схемы Рози} -- граф с теми
же вершинами и рёбрами, рёбра графа пронумерованы, а слов на них
нет. {\it Протокол детерминированной эволюции} --
последовательность облегчённых нумерованных схем Рози, номеров
соответствующих опорных рёбер, а также информация о минимальных
путях, содержащих опорное рёбро и не являющихся допустимыми.
\end{definition}

\section{Доказательство основной теоремы}

Напомним, что {\em пословнная сложность}  $P(N)$ сверхслова $W$ --
количество различных подслов $W$ длины $N$. {\em Показатель
рекуррентности} $P_2(N)$ для равномерно рекуррентного сверхслова
$W$ -- минимальное число $P_2(N)$ такое, что в любом подслове
сверхслова $W$ длины $P_2(N)$ встретятся все подслова $W$ длины
$N$.

\begin{lemma} \label{Lm6_4}
а) Если $W$ --  морфическое слово, порождённое примитивным
морфизмом $\varphi$, то его показатель рекуррентности $P_2(N)$
ограничен сверху  линейной функцией: $P_2(N)\leq C_3 N$.

б) Если $W$ -- сверхслово с не более, чем линейным показателем
рекуррентности, то пословная сложность $W$ также не более чем
линейна.
\end{lemma}

Нам потребуется следующий результат Ж.~Кассиня

\begin{lemma}[\cite{Cassaigne}] \label{compl}
Если для р.р. сверхслова $W$ существует такая константа $C$, что
для всех $N$ выполнено $P(N)\leq C N$, то функция $P(N+1)-P(N)$
ограничена.
\end{lemma}

\begin{definition}
{\em Масштаб схемы} -- наименьшая из длин слов опорных рёбер.
\end{definition}

Из лемм \ref{Lm6_4} и \ref{compl} вытекает

\begin{lemma} \label{cas}
а) Если сверхслово $W$  морфично, то количество вершин суммарной
степени более $2$ во всех графах Рози $G_k(W)$ ограничено.

б) Количество вершин в любой схеме Рози для $W$ ограниченно.
Количество различных облегчённых нумерованных схем, возникающих
при детерменированной эволюции, конечно. Существует такая
константа $E$, что для любой схемы $S$ длины всех слов на рёбрах
этой схемы не превосходят $EM$, где $M$ -- масштаб схемы.
\end{lemma}

\begin{definition}
Пусть $A\sqsubseteq W$, $S$ -- схема Рози для $W$. Множество
симметричных путей в $S$, слова которых являются подсловами $A$,
обозначим $S(A)$. Если $S(A)$ непусто, в этом множестве есть
максимальный элемент относительно сравнения $\sqsubseteq $,
который будем называть {\em нерасширяемым путём}. Очевидно, для
каждого пути из $S(A)$ есть путь, который содержит его и является
нерасширяемым. Назовём этот путь {\em максимальным расширением}.
{\em Набор проверочных слов порядка $k$} -- это набор слов
$\{q_i\}$, в который входят $\psi(\varphi^k(a_i))$ для всех букв
алфавита $\{a_i\}$, а также $\psi(\varphi^k(a_ia_j))$ для
всевозможных пар последовательных букв слова
$\varphi^{\infty}(a_1)$.
\end{definition}

\begin{proposition} \label{pr7_5}
Существуют такое  $\lambda_0>1$, что в наборе проверочных слов с
номером $k$ длина любого слова $q$ удовлетворяет двойному
неравенству $D_{1}\lambda_0^k<|q|<D_{2}\lambda_0^k$ для некоторых
$D_1$ и $D_2$.
\end{proposition}

\begin{definition}
{\em Оснастка} $(S,k)$ порядка $k$ определяется для
пронумерованной схемы Рози $S$. Для получения оснастки берётся
набор проверочных слов порядка $k$, для каждого проверочного слова
$v_i$ рассматривается набор $S(q_i)$. Каждый путь в схеме задаётся
упорядоченным набором чисел -- номеров рёбер схемы. Таким образом,
$S(q_i)$ задаётся множеством таких упорядоченных наборов, а
оснастка получается, если взять такие наборы для всех проверочных
слов и облегчённую нумерованную схему $S$. {\em Размер оснастки}
-- максимальная длина (в рёбрах) по всем путям из $S(q_i)$ для
всех $q_i$ -- проверочных слов порядка $k$.
\end{definition}

\begin{lemma} \label{size}
а) Размер оснастки отличается от масштаба схемы в ограниченное
число раз. Количества типов оснасток ограниченно. И, кроме того,
оснастка $(S,k+1)$ является функцией оснастки $(S,k)$.

б) Если размер оснастки $(S,k)$ достаточно большой, то по оснастке
$(S,k)$ однозначно определяется оснастка $(\Evol(S),k)$.
\end{lemma}

{\bf Теорема \ref{ThMain}} для случая морфического  $W$ с
примитивным порождающим морфизмом вытекает из ограничености типов
оснасток и однозначности перехода к следующей оснастке, общий
случай требует привлечения результатов Ю.~Л.~Притыкина
\cite{SemenovReview} а также дополнительного соображения:

\begin{lemma}[\cite{Ehren}] \label{eh}
Следующие три условия эквивалентны: а)  Слово $\varphi
^{\infty}(a)$ не является равномерно рекуррентным. б) В $\varphi
^{\infty}(a)$ есть бесконечно много $\varphi$-ограниченных
подслов. в) Существует непустое $w\in A^*$ такое, что $w^n$
является подсловом $\varphi ^{\infty}(a)$ для любого $n$.
\end{lemma}


\begin{thebibliography}{aaa}

\bibitem{SemenovReview}
 {\sl
Ан.~А.~Мучник, Ю.~Л.~Притыкин, А.~Л.~Семенов.} {\it
Последовательности, близкие к периодическим.} УМН, 64:5(389)
(2009), 21-96


\bibitem{Ehren} {\sl A. Ehrenfeucht and G. Rozenberg.} {\it Repetition of subwords in
$DOL$ languages.} Information and Control, 59(1--3):13--35, 1983.

\bibitem{Frid} {\sl  А. Э. Фрид}, {\it О графах подслов
$DOL$-последовательностей.} Дискретн. анализ и исслед. опер., сер.
1, 6:4 (1999), 92--103


\bibitem{BelKondMitr} {\sl A.Ya.Belov, G.V.Kondakov, I.Mitrofanov.}
{\it Inverse problems of symbolic dynamic}.

\bibitem{BelovChernInt}
{\sl A.Ya.~Kanel-Belov, A.L.~Chernyat'ev.} {\it Describing the set
of words generated by interval exchange transformation}. Comm. in
Algebra, Vol. 38, No 7, July 2010, pages 2588--2605.

\bibitem{Cassaigne} {\sl J.Cassaigne.} {\it Special factors with linear
subword complexity.} Developments in language theory, II
(Magdeburg, 1995), 25-34, World Sci. Publ., River Edge, NJ, 1996.


\bibitem{FZ3} {\sl Ferenczi, L. Zamboni}, {\it Languages of $k$-interval exchange
transformations.} http://iml.univ-mrs.fr/~ferenczi/fz3.pdf
Bulletin London Math. Soc. 40 (2008), p. 705--714.

\bibitem{Ferenci2} {\sl S.Ferenczi} {\it
Dynamical generalizations of the lagrange spectrum.}
http://iml.univ-mrs.fr/~ferenczi/lagiet.pdf.




\bibitem{Ra3} G. Rauzy, \'Echanges
d'intervalles et transformations induites, (in French), Acta
Arith. 34 (1979), p. 315-328.


\bibitem {VL}
 Vershik, A. M.; Livshits, A. N.
{\it Adic models of ergodic transformations, spectral theory,
substitutions, and related topics. Representation theory and
dynamical systems}, 185--204, Adv. Soviet Math., 9, Amer. Math.
Soc., Providence, RI, 1992.


\end{thebibliography}
\end{document}